\documentclass[reqno]{amsart}
\usepackage{newcent}       
\usepackage{helvet}         
\usepackage{courier}        
\usepackage{amsfonts}
\usepackage{amsfonts,amssymb,amsmath}
\usepackage[latin1]{inputenc}
\usepackage{color}
\usepackage{graphicx}
\usepackage{amsmath}

\usepackage{mathrsfs}
\usepackage{mathpazo}
\usepackage{tikz}

\makeatletter
\@addtoreset{equation}{section}
\makeatother

\renewcommand\thefigure{\thesection.\@arabic\c@figure}
\renewcommand\thetable{\thesection.\@arabic\c@table}

\newcommand{\mc}[1]{{\mathcal #1}}

\newcommand{\bb}[1]{{\mathbb #1}}

\newcommand{\<}{\langle}
\renewcommand{\>}{\rangle}

\definecolor{DBlue}{rgb}{0.1,0,0.55}
\definecolor{DRed}{rgb}{0.55,0,0.1}

\def\emptysquare{{\hbox{\vrule height6pt width0.6pt depth0pt%
\vbox{\hrule height0.6pt width4.8pt depth0pt%
\vglue4.8pt%
\hrule height0.6pt width4.8pt depth0pt}%
\vrule height6pt width0.6pt depth0pt}}}

\def\qed{\unskip\nobreak
\hfil\penalty50\hskip1.75em\null\nobreak\hfil\emptysquare
{\parfillskip=0pt \finalhyphendemerits=0 \par}\medskip}

\author{Patr\'{\i}cia Gon\c{c}alves}
\address{Departamento de Matem\'atica, PUC-RIO, Rua Marqu\^es de S\~ao Vicente, no. 225, 22453-900, Rio de Janeiro, Rj-Brazil and CMAT, Centro de Matem\'atica da Universidade do Minho, Campus de Gualtar, 4710-057 Braga, Portugal.}
\email{patricia@mat.puc-rio.br}
\author{Milton Jara}
\address{IMPA, Estrada Dona Castorina 110, Jardim Bot\^anico, CEP 22460-340, Rio de Janeiro, Brazil}
\email{mjara@impa.br}
\title[The Einstein relation for the KPZ equation]{The Einstein relation for the KPZ equation}
\date{\today}

\begin{document}

\begin{abstract}
We compute the non-universal constants in the KPZ equation in one dimension, in terms of the thermodynamical quantities associated to the  underlying microscopic dynamics. In particular, we derive the second-order Einstein relation $\lambda = \frac{1}{2}\frac{d^2}{d\rho^2} \chi(\rho) D(\rho)$ for the transport coefficient $\lambda$ of the KPZ equation, in terms of the conserved quantity $\rho$, the diffusion coefficient $D$ and the static compressibility of the system $\chi$.
\end{abstract}
\maketitle

\section{Introduction}
One of the most challenging problems in statistical mechanics corresponds to the study of the evolution of nonequilibrium systems, and in particular the derivation of effective equations in terms of the relevant thermodynamical quantities of the models in question. One particular problem which has received a lot of attention recently, is the evolution of random growing interfaces governed by {\em local} stochastic rules. In the seminal work \cite{KPZ}, Kardar, Parisi and Zhang proposed an effective equation, nowadays wide known as the KPZ equation, for the evolution of the fluctuations around the mean of a flat growing interface. They argued that the evolution of a fluctuating interface is governed by three competing factors: {\em roughening}, represented by the presence of a space-time white noise, {\em smoothing}, represented by a diffusive term appearing in the form of a Laplacian operator, and a slope-dependent {\em growth}, represented by a nonlinear transport term. Taking these three ingredients into account, \cite{KPZ} proposed the equation
\begin{equation}\label{KPZ}
\partial_t h = D \Delta h +  \lambda (\nabla h)^2 + \sqrt{2D\chi} \mc W,
\end{equation}
where $\mc W$ is a normalized space-time white noise, that is, a Gaussian noise with correlations given by $\<\mc W(x,t) \mc W(x',t')\> = \delta(x,x')\delta(t,t')$ and $D,\lambda,\chi$ are constants.\\

The main result we want to report here is the computation of the constant $\lambda$ in terms of the constants $D$ and $\chi$. We will be more specific later on, here we just mention that $D$ is the {\em diffusivity} of the system, and $\chi$ is the {\em static compressibility} of the system. By scaling arguments, the quadratic dependence on the slope $\nabla h$ is the only nontrivial quantity which is compatible with both the Laplacian term and the white noise.\\

In \cite{KPZ}, starting from the KPZ equation the authors argued that a growing interface has a statistically self-similar structure with universal scaling exponents for its width and for its spatial correlations. In particular, they predicted, starting from the KPZ equation above, that a one-dimensional interface has fluctuations of order $t^{1/3}$, in contrast with the Edwards-Wilkinson (EW) exponent $t^{1/4}$ of the fluctuations of interfaces in equilibrium. In the EW class the fluctuations evolve according to the Ornstein-Uhlenbeck equation, which corresponds to taking $\lambda=0$ in the equation \eqref{KPZ}. In principle, the aforementioned feature should be shared by any discrete or continuous, one-dimensional  growth models with local stochastic dynamics. Starting from \cite{vBKS}, this question has been investigated through Monte Carlo simulations for various simplified models, like the Eden model, the random deposition model, and the polynuclear growth model. Universal exponents were confirmed by those simulations, giving support to the KPZ conjecture. We say that a model belongs to the KPZ universality class if its corresponding fluctuations follow the exponents predicted by \cite{KPZ}. For an early review and further references, we refer to \cite{Kru}.\\

A first ground-breaking contribution, and also the first mathematically rigorous result in this direction, was obtained by \cite{BG}. In that article, the authors derived the so-called Cole-Hopf solution of the KPZ equation as the scaling limit of the fluctuations of an interface associated to  a particular interacting particle system: the {\em weakly asymmetric} simple exclusion process. The simple exclusion process is a system of interacting particles evolving on $\bb Z$ with the following stochastic dynamics. Let $p,q \geq 0$ be such that $p+q=1$. Independently, the position of a particle is updated with rate 1. The update rules are the following: the particle tries to jump to the right neighboring site with probability $p$ and to the left neighboring site
with probability $q$. The jump is successful if there is no particle at the target position. The weakly asymmetric scaling corresponds to the choice $p=\frac{1}{2}+a\sqrt \epsilon,$ where $\epsilon >0$ is a scaling parameter which goes to 0 and $a$ can be understood as the strength of the asymmetry.\\

The work of \cite{BG} does not say anything about the scaling exponents of the fluctuations of the system, but it clarifies the role of the KPZ equation as an effective model for the evolution of fluctuations in one-dimensional growth models. In particular, the behavior of asymmetric growth models like the TASEP (which corresponds to the choice $p=1$ above), the polynuclear growth model, etc, should be related with the {\em long-time behavior} of the KPZ equation.\\

A second ground-breaking contribution is contained in the paper \cite{Joh}, where it is computed the probability distribution of the height function on a discrete growth model, known as the single step, with a particular initial configuration known as the wedge profile (see also \cite{BDJ,PS}).  There it is  shown that the fluctuations of the height function of that
 model are given by the Tracy-Widom (TW) distribution \cite{TW}, and that the corresponding scaling exponent is effectively $\frac{1}{3}$. We point out that the TW distribution is observed in a strongly asymmetric regime, which in the case of exclusion processes corresponds to the choice  $p\neq \frac{1}{2}$ independently of the scaling. In \cite{ACQ,SS}, the authors showed that the KPZ equation serves as a crossover equation connecting the KPZ universality class ($p \neq \frac{1}{2}$ in the case of the exclusion process) and the EW universality class ($p=\frac{1}{2}$ in the case of the exclusion process). More precisely, they showed that, as $t \to \infty$, $t^{-1/3} \bar{h}(t,0)$ converges to a TW distribution, and as $t \to 0$, $t^{-1/4} \bar{h}(t,0)$ converges to a normal distribution. Here $\bar{h}(t,0)$ is the centered solution of the KPZ equation with wedge initial profile. In this sense, the KPZ equation is the universal object that serves as separation between the EW ($t^{1/4}$-scaling) and KPZ ($t^{1/3}$-scaling) universality classes.\\

 More recently, there has been a new way to solve the KPZ equation. The analysis is performed model-by-model and the idea is to obtain some determinantal formulas related to some functionals of the model in question. These formulas can be solved by using the machinery of random matrix theory, which provides very detailed information on the limiting distribution of those functionals \cite{Co}. \\

The main drawback of all the approaches described above is the lack of generality. All of them are based on intricate combinatoric properties of the model considered; in \cite{Spo} the term {\em stochastic integrability} is coined to point out this fact. In particular, the non-universal constants on the KPZ equation do not depend on the thermodynamical properties of the proposed models. This is not true in general, as we will see below. From the point of view of the thermodynamical properties of the microscopic growth models, it is more convenient to adopt a different approach. For that purpose, define $\mc Y_t = -\nabla h_t.$ Since $h_t$ solves the KPZ equation \eqref{KPZ}, then the field $\mc Y_t$ solves the {\em stochastic Burgers equation} (SBE)
\begin{equation}\label{SBE}
\tag{SBE}
\partial_t \mc Y_t = D \Delta \mc Y_t + \lambda \nabla \mc Y_t^2 + \sqrt{2\chi D} \nabla \mc W.
\end{equation}
If we start with a discrete growth model $\{\zeta_t(x); x \in \bb Z\}$ with local interactions, the discrete gradients defined as $\eta_t(x) = \zeta_t(x) - \zeta_t(x+1)$ can be interpreted as a {\em conservative} interacting particle system, and vice-versa, for a given conservative interacting particle system, the cumulative currents of the system can be interpreted as an interface growth model. Therefore, from now on we stick to the interacting particle systems interpretation.\\

An important feature which is not captured by the results of \cite{BG}\footnote{
We point out that in \cite{KPZ}, the authors use generic constants in front of the three terms of this equation, and they do not discuss their meanings in terms of thermodynamical quantities of the underlying systems.}, is the dependence of $\lambda$ in the thermodynamical
 properties of the system in consideration. At principle there is no obvious relation between $\lambda$ and the underlying particle system, but in fact, there is a relation between them and the purpose of this paper is to describe precisely that relation. Notice, that the steady states of the exclusion process are parametrized by its density $\rho$. Taking the jump rates $p=\frac{1}{2}+a\sqrt \epsilon$, we get to the SBE equation above with $\lambda = a$, which does not depend on $\rho$. This is very particular of the chosen jump rates and in general is not true. We will show here a {\em second-order Einstein relation}, namely the relation
$$\frac{d\lambda}{da} = \frac{1}{2} \frac{d^2}{d\rho^2} {\chi D}.$$

The outline of this paper is as follows. In section \ref{sec model} we present the model that we have chosen in order to obtain our main result, namely,  the gradient Kawasaki dynamics. In Section \ref{sec kpz} we describe its equilibrium fluctuations for the KPZ scaling. In section \ref{conclusion} we present some conclusions and future directions.

\section{The gradient Kawasaki dynamics}\label{sec model}

Consider the space $\Omega = \{0,1\}^{\bb Z}$ of binary sequences. We call $\eta = \{\eta(x); x \in \bb Z\}$ the elements of $\Omega$ and we interpret $\eta(x)=1$ as having a particle at site $x$ and $\eta(x)=0$ as having a hole at the site $x$. Given a Gibbs measure $\mu$ in $\Omega$, the Kawasaki dynamics is a local, particle-conservative interacting particle system for which the measure $\mu$ is invariant and reversible. Given a Hamiltonian $\mc H$, in \cite{Nag} it is explained how to choose a Kawasaki dynamics satisfying an additional property, the so-called {\em gradient condition}. Let $J = \{J_A; A \subseteq \bb Z\}$ be a finite-range, translation-invariant potential, that is, there exists a constant $R$ such that $J_A=0$ whenever diam$(A) > R$ and $J_A=J_{A+x}$ for any $x \in \bb Z$. For $\Lambda\subset {\mathbb{Z}}$, define the Hamiltonian $\mc H_\Lambda$ as $$\mc H_\Lambda(\eta) = \sum_{A \cap \Lambda \neq \emptyset} J_A \eta(A),$$ where $\eta(A) = \prod_{x \in A} \eta(x)$. Consider an inverse temperature $\beta$ such that the Gibbs measure $$\mu_\beta = \lim_{\Lambda} \frac{1}{Z_{\Lambda,\beta}} e^{-\beta \mc H_\Lambda}$$ is well-defined. Above, $Z_{\Lambda,\beta}$ is a normalizing constant. Let us denote by $\eta^{0,1}$ the configuration obtained from $\eta$ by exchanging the occupation numbers of $0$ and $1$, namely
$$\eta^{x,x+1}(y)=\eta(x+1)\textbf{1}_{y=x}+\eta(x)\textbf{1}_{y=x+1}+\eta(x)\textbf{1}_{y\neq x,x+1}.$$
Let $\theta$ be the standard spatial shift in $\bb Z$, that is, for $x\in\mathbb{Z}$ and $\eta\in\Omega$, $$\theta_x(\eta)(y)=\eta(x+y).$$ In \cite{Nag}, it is shown the existence of a function $c: \Omega \to \bb R$ such that:
\begin{itemize}
\item[i)]{\em Local dynamics.} The value of $c(\eta)$ depends only on $\{\eta(x); |x| \leq r\}$ for some $r >0$,

\quad

\item[ii)] {\em Ergodicity and exclusion rule.} $c(\eta)>0$ if $\eta(0) \neq \eta(1)$ and $c(\eta)=0$ if $\eta(0)=\eta(1)$,

\quad

\item[iii)] {\em Detailed balance.} If $\eta(0) \neq \eta(1)$,
\[
\frac{c(\eta)}{c(\eta^{0,1})} = e^{-\beta(\mc H(\eta^{0,1})- \mc H(\eta))},
\]

\item[iv)] {\em Gradient condition.}  There exists a local function $\omega: \Omega \to \bb R$ such that
$$c(\eta)(\eta(0)-\eta(1)) = \omega(\eta)-\omega(\theta_1 \eta).$$
\end{itemize}
Since $c$ is local, it is bounded. Therefore, without loss of generality, we can assume that $c$ is  bounded by 1. Now we define the Kawasaki dynamics associated to $c$.

Independently for each $x \in \bb Z$, we exchange the occupation variables $\eta_t(x)$ and $\eta_t(x+1)$ with rate 1. The exchange is accepted with probability $c(\theta_{-x} \eta_t)$ and otherwise rejected. We call the stochastic process $\{\eta_t\,: \,t \geq 0\}$ defined in this way the {\em gradient Kawasaki dynamics} associated to $c$. The gradient Kawasaki dynamics is conservative in the sense that particles are nor destroyed neither created by the dynamics. This conservation law implies that a family of invariant measures can be obtained introducing a {\em fugacity}: define $$\mc H_{\Lambda,\beta, \phi}(\eta) = \beta \mc H_\Lambda(\eta) + \phi \sum_{x \in \Lambda} \eta(x).$$ Then, the measure $$\mu_{\beta,\phi} = \lim_{\Lambda} \frac{1}{Z_{\Lambda,\beta,\phi}} e^{-\mc H_{\Lambda,\beta,\phi}}$$ is invariant under the Kawasaki dynamics described above.

We can introduce an asymmetry into this dynamics by defining $$c_\gamma(\eta) = c(\eta)(1-\gamma \eta(1)(1-\eta(0)))$$
for $\gamma \in (0,1]$. It turns out that the measures $\mu_{\beta, \phi}$ are also invariant with respect to the perturbed dynamics, that is, the dynamics associated to $c_\gamma$. In fact, in \cite{Nag} it is shown that this property is equivalent to the gradient condition stated above.

Let us define the {\em density of particles} as $$\rho(\phi) = \int \eta(0) \mu_{\beta,\phi}.$$ Since particles are conserved by the dynamics, it is natural to consider the fugacity as a function of the density of particles, $\phi = \phi(\rho)$ and to reparametrize the invariant measures by the density of particles: we write $\nu_\rho = \mu_{\beta,\phi(\rho)}$. The two basic thermodynamic quantities associated to the symmetric dynamics are the {\em diffusivity}, defined as $$ D(\rho) = \frac{d}{d\rho} \int \omega d\nu_\rho,$$ and the {\em static compressibility}, defined as
\[
\chi(\rho) = \lim_{n \to \infty} \int \Big(\frac{1}{\sqrt{n}} \sum_{x=1}^n(\eta(x)-\rho)\Big)^2 d \nu_\rho.
\]

In the asymmetric case, another meaningful quantity is given by the {\em flux function} $$H(\rho) = \int j d \nu_\rho,$$ where $j(\eta) = c(\eta) \eta(0)(1-\eta(1)).$
One of the versions of the fluctuation-dissipation relation states that $$\int c(\eta)(\eta(1)-\eta(0))^2 d\nu_{\rho} = 2 \chi(\rho) D(\rho).$$ Due to the gradient condition, we have the relation
\[
c(\eta)\eta(1)(1-\eta(0)) = \tfrac{1}{2}c(\eta)(\eta(1)-\eta(0))^2 + \tfrac{1}{2} (\omega(\theta_1 \eta ) -\omega(\eta)),
\]
and in particular, $H(\rho) = \chi (\rho)D(\rho)$.

\section{The KPZ scaling}\label{sec kpz}

Let $\epsilon \in (0,1)$ be a scaling parameter which will go to $0$. As pointed out in \cite{BG}, the stochastic Burgers equation appears as the scaling limit for the density of particles on a {\em weakly asymmetric}, conservative particle system, under a {\em diffusive scaling}. This corresponds in our case to the choice $$\gamma = a \sqrt{\epsilon}.$$ Fix $\rho \in (0,1)$. Let $\{\eta_{t\epsilon^{-2}}\, :\, t \geq 0\}$ be the gradient Kawasaki dynamics with initial distribution $\nu_\rho$, associated to $c_\gamma$ for $\gamma =a \sqrt{\epsilon}$. Notice that the system is being speeded up in the diffusive time scaling. Since the solutions of the stochastic Burgers equation are distribution-valued, it is convenient to define the density fluctuation field $\{\mc Y_t^\epsilon; t \geq 0\}$ through its action over test functions $F$ living in the Schwartz space $\mathcal{S}(\mathbb{R})$ as:
\[
\<\mc Y_t^\epsilon, F \> = \sqrt{\epsilon}\sum_{x \in \bb Z}  (\eta_{t\epsilon^{-2}}(x)-\rho)F(\epsilon x).
\]
For $x\in\mathbb{Z}$, $\eta\in\Omega$ and $f: \Omega \to \bb R$ let us introduce the notation $f_x(\eta)  = f(\theta_{-x} \eta)$.
Using the Markovian character of the evolution of $\eta_{t\epsilon^{-2}}$, the total time-derivative of $\<\mc Y_t^\epsilon, F \> $ is equal to the sum of three terms: a diffusive term
\[
\sqrt{\epsilon} \sum_{x \in \bb Z} \omega_x(\eta_{t \epsilon^{-2}}) \Delta F(\epsilon x),
\]
a transport term
\[
a \sum_{x \in \bb Z} j_x(\eta_{t\epsilon^{-2}}) \nabla F(\epsilon x)
\]
and a noise term of instantaneous variance given by
\[
\epsilon \sum_{x \in \bb Z} c_x(\eta_{t \epsilon^{-2}})(\eta_{t \epsilon^{-2}}(x+1)-\eta_{t \epsilon^{-2}}(x))^2 (\nabla F(\epsilon x))^2.
\]
Above, $\Delta F(\epsilon x)$ and $\nabla F(\epsilon x)$ denote the second and the first space derivative of $F$ in $\epsilon x$, respectively.
In order to identify the stochastic partial differential equation ruling the space-time evolution of limit of $\mathcal Y_t^{\epsilon}$, we need to close the terms above as functions of $\mathcal{Y}_t$. For that purpose, we need to use what is called in the literature as the Boltzmann-Gibbs principle, which was first introduced by \cite{BR} and was proved in this context by \cite{D-MPSW}. This principle allows to replace the diffusive term by
\[
D (\rho)\sqrt \epsilon \sum_{x \in \bb Z} (\eta_{t \epsilon^{-2}}(x)-\rho) \Delta F(\epsilon x) = \<\mc Y_t^\epsilon,D(\rho) \Delta F\>.
\]
The ergodic theorem shows that the instantaneous variance of the noise term is well approximated by $2 \chi(\rho) D (\rho)\int (\nabla F)^2 dx$. The main novelty comes from the analysis of the transport term. Notice in first place that there is no $\sqrt \epsilon$ factor in front of the sum. Therefore, at first glance it is not reasonable to think that this term is bounded. It turns out that the Boltzmann-Gibbs principle shows that this term is well approximated by
\[
aH'(\rho) \sum_{x \in \bb Z} (\eta_{t \epsilon^{-2}}(x)-\rho) \nabla F(\epsilon x) = a\epsilon^{-\frac{1}{2}} H'(\rho) \<\mc Y_t^\epsilon, \nabla F\>.
\]
Therefore, we need to impose $H'(\rho)=0$ if we want to see a non-trivial limit of $\mc Y_t^\epsilon$. As already noticed by \cite{KPZ}, this is a natural assumption, since the fluctuations should be observed around the {\em characteristic lines} of the system. Moreover it is not a real restriction, because after a Galilean transformation of the field, the thermodynamical constants $\chi(\rho)$, $D(\rho)$ remain unchanged, while the new flux satisfies $$\tilde H'(\rho) = H'(\rho) -v,$$ where $v$ is the velocity of the Galilean transformation.

Notice that, since we are assuming $H'(\rho)=0$, the usual Boltzmann-Gibbs principle does not give any useful information about the behavior of the transport term. In order to study the limit of the transport term, in \cite{GJ,GJ2, GJ3,GJS} we introduced what we call the {\em second-order Boltzmann-Gibbs principle}, see for example Theorem 7  in \cite{GJ3}. For that purpose, let $\iota: \bb R \to \bb R$ be a positive test function with $\int \iota(x) dx =1$, and define the approximation of the identity $\iota_\delta(x) = \frac{1}{\delta}\iota(\frac{x}{\delta})$, $\delta >0$. Let $f:\Omega \to \bb R$ be a local function and $$\tilde{f}(\rho) = \int f d\nu_\rho.$$ Assume that $\tilde{f}'(\rho) = 0$. Then,  the second-order Boltzmann-Gibbs principle asserts that for any test function $F\in\mathcal{S}(\mathbb{R})$, the average
\[
\sum_{x \in \bb Z} f_x(\eta_{t \epsilon^{-2}}) F(\epsilon x)
\]
is very well approximated by
\[
\tfrac{1}{2}\tilde{f}''(\rho) \epsilon \sum_{x \in \bb Z} F(\epsilon x) (\<\mc Y_t^\epsilon, \iota_\delta^{\epsilon x}\>^2-\frac{\kappa\chi(\rho)}{\delta}),
\]
where we denote by $\iota_\delta^{\epsilon x}$ the approximation of the identity centered at $\epsilon x$ and $\kappa = \int \iota(x)^2 dx$. The approximation holds when $\epsilon \to 0$ {\em and then} $\delta \to 0$. Notice the Wick renormalization represented by the diverging term $\frac{\kappa\chi(\rho)}{\delta}$. The justification of this approximation comes from the equivalence of ensembles. Fix a local function $f:\Omega  \to \bb R$ and let $\psi_f^\ell(\sigma)$ be the expectation of the function $f$ with respect to the canonical ensemble in $\Lambda_\ell = \{1,\dots,\ell\}$ with density of particles $\frac{1}{\ell} \sum_{x=1}^\ell \eta(x) = \sigma$:
\[
\psi_f^\ell(\sigma)=E_\rho^\ell\Big[f\Big|\frac{1}{\ell} \sum_{x=1}^\ell \eta(x) = \sigma\Big],
\]
where $E_\rho^\ell$ denotes the expectation with respect to the canonical ensemble on $\Lambda_\ell$.
 Then, morally, this condition expectation satisfies a "Taylor" expansion, in the sense that
\[
\psi_f^\ell(\sigma) = \tilde{f}(\rho) + \tilde{f}'(\rho)(\sigma-\rho) + \tfrac{1}{2}\tilde{f}''(\rho) (\sigma-\rho)^2
\]
plus an error term of order $(\sigma-\rho)^3$.

Applying the second-order Boltzmann-Gibbs principle for $f = j$, we see that the transport term is well approximated by
\[
\tfrac{a}{2} H''(\rho)\epsilon \sum_{x \in \bb Z} \nabla F(\epsilon x) (\<\mc Y_t^\epsilon, \iota_\delta^{\epsilon x}\>^2-\frac{\kappa\chi(\rho)}{\delta}).
\]
Notice that the Wick renormalisation is not needed in this case, since $\int \nabla F(x) dx =0$. In conclusion, we have just shown that the total time-derivative of $\<\mc Y_t^\epsilon,F\>$ satisfies
\[
\tfrac{d}{dt} \<\mc Y_t^\epsilon,F\>
	= D (\rho)\<\mc Y_t^\epsilon,\Delta F\> + \tfrac{a}{2} H''(\rho)\<(\mc Y_t^\epsilon \ast \iota_\epsilon)^2, \nabla F\>
	+ \sqrt{2\chi(\rho) D(\rho)} \<\mc W_t,\nabla F\>
\]
plus an error term which vanishes as $\epsilon \to 0$ and then $\delta \to 0$. We see that when $\epsilon \to 0$, the process $\mc Y_t^\epsilon$ converges to the process $\mc Y_t$, solution of the stochastic Burgers equation \eqref{SBE} with $\lambda=\frac{a}{2}H''(\rho)$.
Recall that $a$ represents the strength of the asymmetry introduced in the system. In particular, we identify the constant $\lambda$ of the KPZ/stochastic Burgers equation with $$\lambda=\frac{a}{2}H''(\rho),$$ which proves the second-order Einstein relation for the KPZ/stochastic Burgers equation. From all these conclusion, the constants in the KPZ equation are given in terms of the thermodynamical quantities of the system and the SBE equation for $\mathcal{Y}_t$ now reads as

\begin{equation}
\partial_t \mc Y_t = \tilde \omega'(\rho) \Delta \mc Y_t + \frac{a}{2}\tilde j\,''(\rho) \nabla \mc Y_t^2 + \sqrt{2\chi(\rho) \tilde \omega'(\rho)} \nabla \mc W.
\end{equation}

\section{Conclusions and comments}\label{conclusion}

We have shown that the non-universal coefficients of the KPZ equation in $d=1$ can be obtained as functions of the thermodynamical quantities associated to the underlying interacting particle systems. More precisely, denote by $\rho$ the average value of the conserved quantity in a one-dimensional, conservative, weakly asymmetry stochastic dynamics and by $\nu_\rho$ the stationary state associated to $\rho$. Let $\chi=\chi(\rho)$ be the static compressibility of the dynamics, which is simply the self-correlation of the conserved quantity with respect to $\nu_\rho$. Let $D=D(\rho)$ be the diffusion coefficient associated to the symmetrized dynamics, computed using the Green-Kubo formula or any other convenient formula. Let $a$ be a parameter regulating the strength of the asymmetry. Then, the space-time fluctuations of the conserved quantity, with respect to the stationary state $\nu_\rho$ are given by the KPZ/SBE equation
\[
\partial \mc Y_t = D \Delta \mc Y_t + \lambda \nabla \mc Y_t^2 + \sqrt{2\chi D} \nabla \mc W,
\]
where $\mc W$ is a space-time white noise of covariance matrix $\delta(x-x')\delta(t-t')$ and $\lambda= \lambda(a,\rho)$ satisfies
\[
\tfrac{d}{da} \lambda = \frac{1}{2} \tfrac{d^2}{d\rho^2} \chi D.
\]
This relation can be understood as an Einstein relation for the KPZ/SBE, since the transport term turns out to be proportional to the strength of the asymmetry, and the constant of proportionality is given in terms of the thermodynamical quantities associated to the model.

\section{Acknowledgements}

PG thanks FCT for support through  the project PTDC/MAT/109844/2009 and to  CNPq (Brazil) for support through the research project ``Additive functionals of particle systems" 480431/2013-2. PG thanks  CMAT for support by ``FEDER" through the
``Programa Operacional Factores de Competitividade  COMPETE" and by
FCT through the project PEst-C/MAT/UI0013/2011. The authors than the warm hospitality of CMAT (University of Minho) where this  work was finished.


\begin{thebibliography}{99}



\bibitem{ACQ}
G.~Amir, I.~Corwin, and J.~Quastel, Probability distribution of the free energy of the continuum
directed random polymer in 1+1 dimensions, Commun. Pure Appl. Math., \textbf{ 64} no.4, 466--537 (2011).

\bibitem{BDJ}
J.~Baik, P.~Deift and K.~Johansson, On the distribution of the length of the longest increasing
subsequence of random permutations, J. Am. Math. Soc., \textbf{12} no. 4, 1119--1178 (1999).

 \bibitem{BG}
 L. Bertini and G. Giacomin, Stochastic Burgers and KPZ equations from particle systems, Comm. Math. Phys., \textbf{183}  no. 3, 571--607 (1997).


\bibitem{BR}
T.~Brox and H.~Rost, Equilibrium fluctuations of stochastic particle systems: the role of
conserved quantities, Ann. Prob.,  \textbf{ 12}, no. 3, 742--759 (1984).

\bibitem{Co}
I. Corwin, The Kardar-Parisi-Zhang equation and Universality Class, Random matrices: Theory and Applications, vol 1, (2012).

\bibitem{D-MPSW}
A.~De~Masi, E.~Presutti, H.~Spohn and W.~D. Wick, Asymptotic equivalence of fluctuation fields for reversible exclusion processes with speed change, Ann. Prob., \textbf{ 14} no.2, 409--423 (1986).


\bibitem{KPZ}
M.~Kardar, G.~Parisi and Y. C. Zhang, Dynamic Scaling of Growing Interfaces, Phys. Rev.
Lett., \textbf{ 56} no. 9, 889--892,  (1986).

\bibitem{Kru}
H. Krug, H. Spohn, in {\em Solids Far from Equilibrium},  Cambridge Univ. Press (1991).




\bibitem{GJ}
P. Gonçalves and M. Jara: Scaling limits of additive functionals of interacting particle systems, Communications on Pure and Applied Mathematics, \textbf{66} no. 5, 649--677 (2013).

\bibitem{GJ2}
P. Gonçalves and M. Jara: Crossover to the KPZ equation, Annales Henri Poincaré, \textbf{13} no. 4, 813--826 (2012).

\bibitem{GJ3}
 P. Gonçalves, P. and M. Jara: Nonlinear fluctuations of weakly asymmetric interacting particle systems, Archive for Rational Mechanics and Analysis, \textbf{212} no.2, 597--644  (2014).


\bibitem{GJS}
P. Gonçalves, M. Jara and S. Sethuraman: A stochastic Burgers equation from a class of microscopic interactions, to appear in Annals of Probability, Arxiv:1210.0017.

\bibitem{Joh}
K.~Johansson, Shape fluctuations and ranodm matrices, Commun. Math. Phys., \textbf{ 209} no.2, 437--476 (2000).

\bibitem{Nag}
Y. Nagahata, The gradient condition for one-dimensional symmetric exclusion process, J. Stat. Phys., \textbf{ 91}, 587--602 (1998).


\bibitem{PS}
 M. Prahofer, H. Spohn, Exact scaling functions for one-dimensional stationary KPZ growth, J. Stat. Phys.,  \textbf{115} no. 1-2, 255--279 (2004).


\bibitem{SS}
T.~Sasamoto and H.~Spohn, The one-dimensional KPZ equation: an exact solution and its universality, Phys. Rev. Lett.,\textbf{ 104}  no.23 (2010).



\bibitem{Spo}
H.~Spohn, Stochastic integrability and the KPZ equation, International Association of Mathematical Physics News Bulletin, 5--10 (2012).


\bibitem{TW}
C.~Tracy and H.~Widom, Asymptotics in ASEP with step initial condition, Commun. Math. Phys., \textbf{ 290} no. 1, 129--154 (2009).

\bibitem{vBKS}
H. van Beijeren, R. Kutner, H. Spohn, Excess  noise for driven diffusive systems, Phys. Rev. Lett., \textbf{54}, 2026--2029 (1985).


\end{thebibliography}
\end{document}